\documentclass[reqno,11 pt]{amsart}
\usepackage{amsmath,amssymb,latexsym}
\usepackage{hyperref}
\usepackage{amsthm}
\usepackage{color}
\usepackage[all,cmtip]{xy}
\usepackage{adjustbox}
\usepackage{verbatim}

\hypersetup{
    colorlinks=true,
    linktoc=all,
        linkcolor=blue,
}

\newcommand{\Z}{{\mathbb{Z}}}
\newcommand{\Q}{{\mathbb{Q}}}

\theoremstyle{plain}

\newtheorem{theorem}{Theorem}[section]
\newtheorem*{theorem*}{Theorem}

\begin{document}

\title{Fruit Diophantine Equation}
\maketitle

\begin{center}
Dipramit Majumdar \& B. Sury
\end{center}

\vspace{3mm}

\section{Introduction}

\noindent A common problem making the rounds in \href{https://thehighergeometer.wordpress.com/2021/07/27/diophantine-fruit/}{social media} (stated
as a problem about distribution of fruits) is the possibility of
finding positive integer solutions in $x,y,z$ of
\begin{equation}\label{fruitdioph}
y^2 +z^2+5 = x^3+xyz.
\end{equation}

In this short note, we prove that 

\begin{theorem}\label{mainthm}
The equation 
$$y^2 - xyz +z^2= x^3 - 5$$
 has no integer solution. 
\end{theorem}

The problem originated from a question ``What is the smallest unsolved diophantine equation? " asked by user Zidane in MathOverflow \cite{MO}, here the size of a polynomial is defined as follows: for $P(x_1,\dots, x_d) = \sum_{i_1,\dots, i_d} a_{i_1, \dots, i_d} x_1^{i_1}\cdots x_d^{i_d} \in \Z[x_1, \dots, x_d]$, define $|P|(x_1,\dots, x_d) = \sum_{i_1,\dots, i_d} |a_{i_1, \dots, i_d}| x_1^{i_1}\cdots x_d^{i_d}$ and the size of $P$ is defined as $h(P)= |P|(2, \dots,2)$. Following this user Bogdan concentrated on the first part of the problem : `` for what smallest $P$ one does not know if there exists integral solution of $P({\bf x})=0$?"

First non-trivial example of such $P$ was obtained in the case when $h(P)=22$ and that corresponded to the equation $x^2+y^2-z^2=xyz-2$. Will Sawin \cite{MO2} showed that this equation has no non-trivial solution. Next `unknown' diophantine equation was  $y(x^3-y)=z^2+2$ with $h(P)=26$ was solved by user Servaes \cite{MSE}. Finally for $h(P)=29$, Bogdan \cite{MSE}  asked solution for two Diophantine equations $y(x^2+2) = 2zx + 2z^2+1$ (which he himself solved later) and $y^2-xyz+z^2 = x^3-5$. The last equation gained quite a bit of popularity in the social media thanks to a blog `theHigherGeometer', where David Roberts posed this problem as a ` fruit equation' \cite{BLOG}. \\

As a consequence of Theorem \ref{mainthm}, we obtain

\begin{theorem}\label{mainthm2}
For any integer $k$, the elliptic curve  $E_k$ given by the Weierstrass equation $Y^2 - kXY = X^3 -(k^2+5)$ has no integral point.
\end{theorem}

\noindent {\bf Acknoledgements:} I would like to thank Aditya Karnataki for sharing this problem on social media which brought it to our attention. We would also like to thank Aubrey de Grey, who pointed out a typo in the previous version. We would like to thank the blog \href{https://thehighergeometer.wordpress.com/}{theHigherGeometer} for making this problem popular. Last, but not the least, we would like to thank mathoverflow user Bogdan for finding out all these diophantine equations.

\vspace{3mm}

\section{Proof of Theorem \ref{mainthm} and \ref{mainthm2}}

{\bf Proof of Theorem \ref{mainthm}} \\

Fix a prospective $x=k$ that gives an integral solution. Then we are looking for an integral solution of the equation

\begin{equation}\label{eq1}
y^2 - kyz + z^2 = k^3-5.
\end{equation}

We divide the proof in 2 cases. \\

{\bf Case 1: $k$ is even.}\\

Then the equation (\ref{eq1})  becomes
$$\bigg(y - \frac{kz}{2} \bigg)^2 - \bigg(\frac{k^2}{4} - 1 \bigg) z^2 = k^3-5.$$
Rewriting in terms of $d := \frac{k}{2}$ and $u = y - \frac{kz}{2}$, this
becomes the Pell's equation
\begin{equation}\label{eq2}
u^2 - (d^2-1)z^2 = 8d^3-5.
\end{equation}
If $d$ is odd, the left hand side and right hand side of equation (\ref{eq2}) are
respectively congruent modulo $8$ to a perfect square and $-5$; this
is impossible.\\
If $d$ is even, the left and right side  of equation (\ref{eq2}) are congruent modulo $4$ to
$u^2 + z^2$ and $-5$ respectively; once again this is impossible. \\
We conclude that there is no solution of equation (\ref{eq1}) with even $k$.\\

{\bf Case 2: $k$ is odd.}\\

 Going modulo 2 of equation (\ref{eq1}), we obtain
 $$y^2+yz+z^2 =0 \pmod{2}.$$
 This implies $y \equiv z \equiv 0 \pmod{2}$. This implies that left side of equation  (\ref{eq1}) is congruent to $0$ modulo $4$. As a consequence, we obtain $k^3-5 \equiv 0 \pmod{4}$ or equivalently $k \equiv 1 \pmod{4}$. Rewriting the equation (\ref{eq1})  as
$$\bigg(y - \frac{kz}{2} \bigg)^2 - (k^2 - 4) \frac{z^2}{4} = k^3-5$$
 
 and substituting $z=2v$ and $u = y-vk$, the equation (\ref{eq1})  becomes the Pell's equation
 \begin{equation}\label{eq3}
u^2 - (k^2-4)v^2 = k^3-5.
\end{equation}

We need to find integral solutions of the Pell's equation $u^2 - (k^2-4)v^2 = k^3-5$, for integers $k \equiv 1 \pmod{4}$. We divide the proof in the following 3 cases.\\

{\bf Case 2a: $k \equiv 1 \pmod{12}$.}\\
In this case, $3 \mid (k+2)$ and $k^3 \equiv 1 \pmod{3}$. Going modulo 3, the equation (\ref{eq3}) gives us $u^2 \equiv 2 \pmod{3}$. As a consequence, we see there is no solution in this case.\\

{\bf Case 2b: $k \equiv 9 \pmod{12}$.}\\
Note that in this case we have $k-2 \equiv 7 \pmod{12}$, hence $6 \nmid (k-2)$.\\
First observe that if all the prime divisors of $k-2$ are congruent to $\pm 1$ modulo 12, then $k-2$ is congruent to $\pm 1 $ modulo $12$, which contradicts the fact that $k-2 \equiv 7 \pmod{12}$. Thus there exists a prime $p$ such that  $p \equiv 5 \text{ or } 7 \pmod{12}$  and $p$ divides $k-2$. This implies $k \equiv 2 \pmod{p}$, hence $k^3 \equiv 8 \pmod{p}$. Going modulo $p$, the equation (\ref{eq3}) gives us $$u^2 \equiv 3 \pmod{p}.$$ 
By quadratic reciprocity, this equation has no solution as $3$ is a square modulo $p$ if and only if $p \equiv \pm 1 \pmod{12}$.\\

{\bf Case 2c: $k \equiv 5 \pmod{12}$.}\\
In this case, we have $3 \mid (k-2)$.
Going modulo $3$, the equation (\ref{eq3}) gives us 
$$u^2 \equiv 0 \pmod{3}.$$ 
Writing $u=3w$ and $k=12r+5$, the equation (\ref{eq3}) becomes
 \begin{equation}\label{eq4}
3w^2 - (4r+1)(12r+7)v^2 = 2^63^2r^3+ 2^43^25r^2+ 2^235^2r +40.
\end{equation}
Going modulo $3$, the equation (\ref{eq4}) gives us 
$$-(r+1)v^2 \equiv 1 \pmod{3}.$$ 
The above equation implies $r \equiv 1 \pmod{3}$, which in turn implies that $k \equiv 17 \pmod{36}$.\\

Writing $k= 36s+17$, we see that $k-2 = 3(12s+5)$. If all the prime divisors of $12s+5 = \frac{k-2}{3}$ are congruent to $\pm 1$ modulo 12, then $\frac{k-2}{3}$ is congruent to $\pm 1 $ modulo $12$, which contradicts the fact that $\frac{k-2}{3} \equiv 5 \pmod{12}$. 
Thus there exists a prime $p$ such that  $p \equiv 5  \text{ or } 7 \pmod{12}$  and $p$ divides $k-2$. This implies $k \equiv 2 \pmod{p}$. \\
Going modulo $p$, the equation (\ref{eq3}) gives us 
$u^2 \equiv 3 \pmod{p}.$ By quadratic reciprocity, this equation has no solution.
\qed

\vspace{3mm}

\noindent Our main theorem \ref{mainthm}, gives us a family of elliptic curves with no integral points.\\

\noindent {\bf Proof of Theorem \ref{mainthm2}.} 

Let $(a,b)$ be an integral solution of the equation $Y^2 - kXY = X^3 -(k^2+5)$. Then $(a,b,k)$ is an integral solution of equation (\ref{fruitdioph}). This is a contradiction as   equation (\ref{fruitdioph}) has no integral solution by Theorem \ref{mainthm}. \qed

\vspace{3mm}

For $|k| \le 5$, Theorem \ref{mainthm2} can also be verified using LMFDB \cite{lmfdb}. We note that $E_k \cong E_{-k}$ over $\Q$ given by $(x,y) \mapsto (x,-y)$. We also remark that $E_k$ need not have trivial Mordell-Weil group.  For example,  using LMFDB \cite{lmfdb} we see that $E_1$  has Mordell-Weil group isomorphic to $\Z$ generated by $\big( \frac{101}{16}, \frac{821}{64} \big)$, on the other hand $E_0, E_2, E_3, E_4$ have trivial Mordell-Weil group. $E_5$ has minimal Weirstrass equation $y^2+xy=x^3 -13x -13$ obtained via change of variable $(x,y) \mapsto (x+2, 2x-y-1)$. Using LMFDB \cite{lmfdb} we see that the elliptic curve $y^2+xy=x^3-13x-13$ has Mordell-Weil group isomorphic to $\Z$ generated by  $\big( -\frac{71}{64}, \frac{593}{512} \big)$, which gives us a rational point (of infinite order) $\big( \frac{-199}{64}, \frac{-4289}{512} \big)$ for $y^2-5xy=x^3-30$.

\end{document}